%% file: LV.tex
\documentclass{article}
\usepackage{amsmath,amsfonts,epsfig,xcolor,comment}
\usepackage{mathptmx}

\usepackage{amsmath,amsfonts,latexsym,amssymb}
\usepackage{mathrsfs}
\usepackage{verbatim}
\usepackage[utf8]{inputenc}
\usepackage[T1]{fontenc}
\usepackage{ae,aecompl}
\usepackage{braket}

\input{./macros}

\newcommand{\ddt}{{\del_t}}

\renewcommand{\hbar}{\bar{{\mathbb H}}^3}

\newcommand{\CC}{\mathbb C}

\newcommand{\R}{\mathbb R}

\newcommand{\Hp}{{{\mathbb H}^2}}
\newcommand{\Hs}{{{\mathbb H}^3}}
\newcommand{\U}{{\mathbb U}}

\newcommand{\T}{{\mathcal T}}

\renewcommand{\htwo}{{\mathbb H}^2}
\newcommand{\Hn}{{{\mathbb H}^n}}
\def\eproof{$\Box$ \medskip}

\newcommand{\ddx}{{\partial_x}}
\newcommand{\ddy}{{\partial_y}}

\renewcommand{\ddz}{{\frac{\partial}{\partial z}}}

\newcommand{\Id}{{\operatorname{id}}}

\newcommand{\hh}{{\mathfrak h}}

\newcommand{\sltc}{\mathfrak{sl}_2\mathbb C}

\newcommand{\pslt}{\mathsf{PSL}_2(\mathbb C)}
\newcommand{\nnn}{{\bf n}}

\renewcommand{\hh}{{\mathfrak h}}

\begin{document}

\title{Variation of holonomy for projective structures and an application to  drilling hyperbolic 3-manifolds}
 \author{Martin Bridgeman\thanks{Research supported by NSF grants DMS-2005498, DMS-1564410 and the Simons Fellowship 675497} \   Kenneth Bromberg\thanks{Research supported by NSF grant DMS-1906095.}}

\date{\today}

\maketitle

\maketitle
\begin{center} 
{\em Dedicated to a great mathematician and friend on the occasion \\
of his sixtieth birthday:  Francois Labourie. \\
}
\end{center}
\begin{abstract}
We bound the derivative of complex length of a geodesic under  variation of the projective structure on a closed surface in terms of the norm of the Schwarzian in a neighborhood of the geodesic. One application is to cone-manifold deformations of acylindrical hyperbolic 3-manifolds.
\end{abstract}

\section{Introduction}
Let $X$ be a Riemann surface or, equivalently, a hyperbolic surface and $\gamma$ a closed geodesic on $X$. A projective structure $\Sigma$ on $X$ determines a holonomy representation of $\pi_1(X)$. If the holonomy of $\gamma$ is not parabolic then $\cL_\gamma(\Sigma)$ is the {\em complex length} of $\gamma$ which is a complex number whose imaginary part is well defined modulo $2\pi$. We let $P(X)$ be the space of projective structures on $X$ and $P_\gamma(X)$ the subspace where the holonomy of $\gamma$ is not parabolic or the identity. Then $\cL_\gamma$ is a smooth function on $P_\gamma(X)$. The goal of this note is to gain quantitative control of the derivative of $\cL_\gamma$.

We begin by describing a formula for the derivative. We first identify the universal cover $\tilde X$ with the upper half plane $\U$ normalized such that the imaginary axis is a lift of $\gamma$. The projective structure $\Sigma$ determines a developing map $f\colon \U\to \chat$ which we semi-normalize so that $f(e^\ell z) = e^{\cL_\gamma(\Sigma)} f(z)$ where $\ell$ is the length of $\gamma$ in $X$. We define the holomorphic vector  $\nnn = z\ddz$ on $\chat$. Then the pull back $f^* \nnn$ is a holomorphic vector field on $\U$ that descends to a vector field on the annular cover $X_\gamma \to X$ associated to $\gamma$. We also denote this vector field on $X_\gamma$ as $f^*\nnn$. Note that $f$ is only determined up to post-composition with an element of $\pslt$ that fixes $0$ and $\infty$. Since $\nnn$ is invariant under any such element we have that $f^*\nnn$  is well defined.

The tangent space $T_\Sigma P(X)$ is canonically identified with $Q(X)$, the space of holomorphic quadratic differentials on $X$. The pairing of a holomorphic vector field with a holomorphic quadratic differential $\phi$ is a holomorphic 1-form so $f^*\nnn \cdot \phi$ is a holomorphic (and hence closed) 1-form on $X_\gamma$. Our formula for the derivative of $\cL_\gamma$ is
\begin{theorem}\label{derivative}
Given $\Sigma \in P_\gamma(X)$ with normalized developing map $f$ and $\phi \in Q(X) \cong T_\Sigma(P_\gamma(X))$ we have
$$d\cL_\gamma(\phi) = -\int_\gamma f^*\nnn \cdot \phi.$$
\end{theorem}

In the special case when $\Sigma$ is Fuchsian then $\nnn = f^*\nnn$ and, as we will see below, the integral can be computed explicitly. In general, to estimate the the integral we need quantitive control over the difference between $\nnn$ and $f^*\nnn$. This is the content of our next result.

Before stating our estimates we define some norms. If $\Phi$ is a quadratic differential on $X$ then $|\Phi|$ is an area form so its ratio with the hyperbolic area form is a non-negative function. We let $\|\Phi(z)\|$ be this function. It is a natural pointwise norm of $\Phi$ and we let $\|\Phi\|_p$ be the corresponding $L^p$-norms with respect to the hyperbolic area form. We also let $\| \cdot \|$ be the hyperbolic length of vector fields. For all norms we write $\| \cdot \|_\gamma$ to represent the sup norm over the curve $\gamma$.

With these definitions we can now state our estimate.
\begin{theorem}\label{normbound}
Let $\Phi$ be a quadratic differential on $\U$ and assume that an $r$-neighborhood (in hyperbolic metric) of the geodesic $\gamma$ given by the imaginary axis  we have $\|\Phi(z)\| \le K$ with $r< 1/2$ and $K/r < 1/4$. Then there exists a locally univalent map $f\colon \U\to \chat$ such that
\begin{itemize}
\item $\Phi = Sf$;

\item $f(it) \to 0$ as $t\to 0$ and $f(it) \to \infty$ as $t\to \infty$;

\item $\|f^*\nnn - \nnn\|_{\gamma} \le \frac{9K}r$ and $\|f_*\nnn - \nnn\|_{f\circ\gamma} \le \frac{9K}r$
\end{itemize}
where $\| \cdot \|_\gamma$ is the supremum of the hyperbolic length of the vector field along the curve $\gamma$  in $\U$ and $\|\cdot\|_{f\circ\gamma}$ is the supremum of the Euclidean length of the vector field along the curve $f\circ\gamma$ in $\CC^*$.
\end{theorem}
The proof Theorem \ref{normbound} is more complicated then one might expect. It involves the use of {\em Epstein surfaces} and estimates in 3-dimensional hyperbolic space. 

An elementary consequence of Theorem \ref{normbound} is the following approximation for $d\mathcal L_\gamma$.
\begin{theorem}\label{main}
Let $\Sigma$ be a projective structure and $\gamma$ a closed geodesic of length $\ell$ such that $\|\Phi\| \leq K$ in the $r$-neighborhood of $\gamma$ with $r\le 1/2$ and $K/r\leq 1/4$. Then 
$$\left|d\cL_\gamma(\phi) + \int_\gamma \nnn\cdot \phi\right| \le \frac{9K\ell}r \|\phi\|_\gamma.$$

\end{theorem}

{\bf Proof:}
We have by Theorems \ref{derivative}  and \ref{normbound}
\begin{eqnarray*}
\left|d\cL_\gamma(\phi) + \int_\gamma \nnn\cdot \phi\right| 
&\leq& \int_\gamma |(f^*\nnn -\nnn)\cdot \phi|\\
&\leq & \|f^*\nnn- \nnn\|_{\gamma} \cdot\|\phi\|_\gamma  \int_\gamma\|dz\|_{\Hp}\\
&=& \frac{9K\ell}r \|\phi\|_\gamma 
\end{eqnarray*}
where $\|dz\|_{\Hp}$ is hyperbolic line element.
\eproof

\subsubsection*{Application to deformations of hyperbolic manifolds}

Let $N$ be an acylindrical hyperbolizable 3-manifold  with boundary $S = \del N$. Then for any (noded) conformal structure $Y$ in the Weil-Petersson completion $\overline{\T(S)}$ of Teichm\"uller space there is a unique geometrically finite hyperbolic structure  $M_Y$ on $N$ with conformal boundary $Y$. The hyperbolic structure $M_Y$ determines a projective structure on $Y$ with Schwarzian quadratic differential denoted $\Phi_Y$. One question that naturally arises is if the $L^2$-norm of $\Phi_Y$ is small, does this imply that the $L^\infty$-norm is also small and therefore the manifold $M_Y$ has almost geodesic convex core boundary. This is not the case as the $L^2$-norm may be small but there are short curves where the $L^\infty$-norm is large. In \cite{wp-paper} we analysed this problem. We showed that for $Y\in \T(S)$ if the $L^2$-norm $\|\Phi_Y\|_2$ is sufficiently small then there is a {\em nearby} $\hat Y \in \overline{\T(S)}$ (in the Weil-Petersson metric) such that the $L^\infty$-norm $\|\Phi_{\hat Y}\|_\infty$ is {\em small}. This $\hat Y$ will generally be a noded surface pinched along curves where the $L^\infty$-norm in $Y$ was not small.

The bounds  in \cite{wp-paper} for  {\em nearby} and {\em small}  were linear in  $\|\Phi_Y\|_2^{\frac1{2n(S) + 3}}$ where $n(S)$ is the maximum number of disjoint geodesics,  in particular for $S$ closed connected of genus $g$ then $n(S) =3g-3$. Our application   is to use  our variation bound in Theorem \ref{main} to significantly improve this estimate and replace the power $\frac1{2n(S) + 3}$ with the constant power $\frac{1}{2}$. Thus the bounds become independent of the topology. Namely we  prove:
\begin{theorem}\label{main_drill}
There exists $K,C > 0$ such that the following holds; Let $N$ be an acylindrical hyperbolizable manifold with boundary $S = \partial N$ and $Y\in \T(S)$ with $\|\Phi_Y\|_2 \leq K$.  Then there exists a $\hat Y \in \overline{\T(S)}$ with 
\begin{enumerate}
\item $d_{\rm WP}(Y, \hat Y) \le C\sqrt{\|\Phi_Y\|_2};$

\item $\|\Phi_{\hat Y}\|_\infty \le C \sqrt{\|\Phi_Y\|_2}.$
\end{enumerate}
\end{theorem}
Our proof  also  holds for {\em relative acylindrical} 3-manifolds as in \cite{wp-paper} but for simplicity we will restrict to the acylindrical case.

\section{Variation of holonomy}
\subsubsection*{Preliminaries}
We consider $X$ a Riemann surface structure on a surface $S$ and $P(X)$ the space of projective structures on $X$. Associated to $\Sigma \in P(X)$ is the holonomy representation  $\rho \in \Hom(\pi_1(S),\pslt)$ (unique up to conjugacy).
The space $P(X)$ can be parametrized by $Q(X)$ the space of holomorphic quadratic differentials by taking the Schwarzian derivative of its developing map $f:\U\rightarrow \hat\CC$. As $Q(X)$ is a vector space, it is its own tangent space, and a variation of $\Sigma$ is given by a quadratic differential $\phi \in Q(X)$. Throughout the paper we will use $\phi$ to denote deformations of projective structures and $\Phi$ (or $\Phi_t$) to denote the Schwarzian of projective structures.

The advantage  of the Schwarzian $\Phi$ is that is uniquely determined by $\Sigma$ while the holonomy representation $\rho$ and developing map $f$ are not. Given a smooth 1-parameter family of projective structures $\Sigma_t$ we get a smooth family $\Phi_t \in Q(X)$ of Schwarzians. While the developing maps $f_t\colon \U\to \chat$ are not uniquely determined we can choose them to vary smoothly and it will be convenient to do this after making some normalizations.

We will be interested in the {\em complex length} of an element $\gamma \in \pi_1(S)$ that represents a closed geodesic for the hyperbolic structure on $X$. After fixing $\gamma$, we identify $\U$ with $\tilde X$ so that the deck action of $\gamma$ on $\U$ is given by $\gamma(z) = e^\ell z$ where $\ell$ is the length of the geodesic representative of $\gamma$. Under the covering map $\U\to X$ the imaginary axis is taken to this geodesic. Our second normalization is to choose the developing map $f$ and corresponding holonomy representation $\rho$ such that $\rho(\gamma)(z) = e^\cL z$. Then $\cL$ is the complex length $\cL_\gamma(\Sigma)$. Note that this is only well defined modulo $2\pi i$.

Extending these normalizations to the 1-parameter family $\Sigma_t$ we get a smooth family of developing maps $f_t$ and holonomy representations with $\rho(\gamma)(z) = e^{\cL_t} z$. Note that while $\cL$ is only well defined modulo $2\pi i$, after this choice is made there is a unique choice of $\cL_t$ that makes the path continuous. In particular the time zero derivative $\dot\cL$ of $\cL_t$ is well defined.

We also let $v$ be the vector field on $\U$ such that $f_*v(z)$ is the tangent vector of the path $f_t(z)$ at $t=0$. The vector field $v$ is not equivariant under the deck action but a standard computation gives
$$v - \beta_* v = f^* \dot\rho(\beta)$$
for all $\beta\in \pi_1(S)$ where $\dot\rho(\beta) \in \sltc$ is the time zero of the $\rho_t(\beta)$. As $\rho_t(\gamma)(z) = e^{\cL_t} z$ we have that $\dot\rho(\gamma) = \dot\cL z \ddz$ so
\begin{equation}
v - \gamma_* v = f^*\left(\dot\cL z\ddz\right).
\label{vector-v}
\end{equation}

If $U$ is an open neighborhood in $\U$ where $f_t$ is injective then $(U, f_t)$ is a projective chart for $\Sigma_t$. In the chart $(U,f)$ the vector field $v$ is represented as $g \ddz$ where $g$ is a holomorphic function. In this chart the {\em Schwarzian derivative} of $v$ is $g_{zzz} dz^2$. This is a quadratic differential  and a computation gives that it is $\phi$ the time zero derivative of the path $\Phi_t$.

\subsubsection*{Model deformations}
We are now ready to prove Theorem \ref{derivative}. Here's the outline:
\begin{itemize}
\item We first construct a model deformation $\phi_\lambda$ on the projective structure on the annulus $X_\gamma$.

\item The model deformation will be very explicit so that we can directly calculate the line integral of $f^*\nnn\cdot\phi_\lambda$ over $\gamma$.

\item We then find a specific $\lambda$ so that the 1-form $f^*\nnn\cdot \phi - f^*\nnn\cdot \phi_\lambda$ is exact. Then the line integrals over both 1-forms will be equal so our previous calculation will give the theorem.
\end{itemize}
A minor issue with this outline is that the model deformation will only be defined on a sub-annulus of $X_\gamma$. We begin with this difficulty.

The set $f^{-1}(\{0,\infty\})$ will be discrete in $\chat$ and $\gamma$-invariant so it will descend to a discrete set in $X_\gamma$. Let $\gamma'$ be a smooth curve, homotopic to $\gamma$, that misses this set and let $A$ be an annular neighborhood of $\gamma'$ that is also disjoint from this set. Then $\tilde A$, the pre-image of $A$ in $\U$, we also be disjoint from $f^{-1}(\{0,\infty\})$. As $\tilde A$ is simply connected we can choose a well defined function $\log f$ on $\tilde A$ and define the vector field $v_\lambda$ on $\tilde A$ by $f_* v(z) = \lambda f(z) \log f(z) \ddz$. We let $\phi_\lambda$ be the Schwarzian of $v_\lambda$ and differentiating $\lambda z \log z$ three times we see that $\phi_\lambda = f^*\left(-\frac\lambda{z^2} dz^2\right)$. Then $\phi_\lambda$ is our model deformation.

Note the choice of $\log$ defines $\cL$ as a complex number rather than just a number modulo $2\pi$. We will use this in the rest of the proof. Next we compute the line integral of $f^*\nnn\cdot \phi_\lambda$ over $\gamma'$:
\begin{lemma}\label{line_calc}
$$\int_{\gamma'} f^*\nnn\cdot \phi_\lambda=-\lambda\cL$$
\end{lemma}

{\bf Proof:} Let $\sigma\colon[0,1] \to \U$ be a smooth path that projects to $\gamma'$ in $X_\gamma$. In particular $\sigma(1) = \gamma(\sigma(0))$ and $f(\sigma(1)) = f(\gamma(\sigma(0)) = e^{\cL} f(\sigma(0))$. Then
\begin{eqnarray*}
\int_{\gamma'} f^*\nnn\cdot \phi_\lambda & = & \int_\sigma f^*\left(z\ddz \cdot -\frac\lambda{z^2} dz^2\right)\\
&=& -\lambda\int_{f\circ \sigma} \frac{dz}{z}\\
&=&- \lambda(\log f(\sigma(1)) - \log f(\sigma(0))\\
&=& -\lambda\cL.
\end{eqnarray*}
\eproof

Next we give a criteria for the form $f^*\nnn\cdot \phi$ to be exact.
\begin{lemma}\label{exact}
Let $v$ be a $\gamma$-equivariant holomorphic vector field on\ \  $\tilde A$ with Schwarzian $\phi$. Then $f^*\nnn\cdot \phi$ is exact on $A$.
\end{lemma}

{\bf Proof:} The vector field $v$ is a section of the holomorphic tangent bundle $T_{\CC}X_\gamma$ so $v_z$ is a function, $v_{zz}$ is a 1-form, $v_{zzz} = \phi$ is a quadratic differential. Therefore  $f^*\nnn \cdot v_{zz} - v_z$ is a holomorphic function on $X_\gamma$ and we'll show that it is a primitive of $f^*\nnn\cdot \phi$.

We can do this calculation in a chart.
Namely choose an open neighborhood $U$ in $\U$ such $f$ is injective on $U$. Then as above $(U,f)$ is a chart for the projective structure and in this chart $v$ has the form $g \ddz$ where $g$ is a holomorphic function on $f(U)$ and $\phi$ is $g_{zzz} dz^2$.  It follows that on this chart $v_z$ is $g_z$, $v_{zz}$ is $g_{zz} dz$ and $\phi = v_{zzz}$ is $g_{zzz} dz^2$. As the derivative of $zg_{zz}- g_z$ is $g_{zzz}$ this shows that $f^*\nnn \cdot v_{zz} - v_z$ is a primitive for $f^*\nnn\cdot \phi$ as claimed. \eproof

{\bf Proof of Theorem \ref{derivative}:} Since $f^*\nnn\cdot \phi$ is holomorphic it is closed and we have
$$\int_\gamma f^*\nnn\cdot \phi = \int_{\gamma'} f^*\nnn\cdot \phi.$$

Next we calculate to see that
$$v_\lambda - \gamma_* v_\lambda = f^*\left(\lambda\cL z\ddz\right)$$
so if $\lambda = \dot\cL/\cL$ then by  \eqref{vector-v} we have that $v - v_{\dot\cL/\cL}$ is $\gamma$-invariant. Therefore by Lemma \ref{exact} the 1-form $f^*\nnn \cdot (\phi - \phi_{\dot\cL/\cL})$ is exact and
$$\int_{\gamma'} f^*\nnn\cdot \phi = \int_{\gamma'} f^*\nnn\cdot \phi_{\dot\cL/\cL} = -\dot\cL$$
where the last equality comes from Lemma \ref{line_calc}. Combining the first and last equality gives the theorem. \eproof

We conclude this section with a simple formula for the  line integral when the projective structure is Fuchsian. For this, rather than represent the annulus as a quotient of the upper half plane $\U$ for this calculation it is convenient to represent the annulus as an explicit subset of $\CC$. Namely, let
$$A_\ell = \left\{z\in\CC \ |\ e^{\frac{-\pi^2}{\ell}}< |z| < e^{\frac{\pi^2}{\ell}}\right\}.$$
This annulus is conformally equivalent $X_\gamma$ when $\ell = \ell_\gamma(X)$ and the circle $|z|=1$ is the closed geodesic of length $\ell$ in $A_\ell$. For this representation of the annulus the vector field $\nnn$ is written as $\nnn = \frac{2\pi i}{\ell}z\ddz$ and we observe that $\nnn$ extends to a holomorphic vector field on all of $\chat$ that is $0$ at $z=0$ and $z=\infty$. 

To decompose the quadratic differential $\phi$ on $A_\ell$ we let $D^+$ and $D^-$ be the disks in $\chat$ with $|z| > e^{\frac{-\pi^2}{\ell}}$ and $|z| < e^{\frac{\pi^2}{\ell}}$, respectively. We recall that any holomorphic function $\psi$ on $A_\ell$ can be written as $\psi_+ +\psi_0 + \psi_-$ where $\psi_+$ and $\psi_-$ extend to holomorphic functions on $D^+$ and $D^-$ that are zero at $z =\infty$ and $z=0$, respectively, and $\psi_0$ is constant. Then $\phi$ can be written as  $\phi = \frac{\psi}{z^2}dz^2$ so the decomposition of $\psi$ gives
\begin{equation}\label{annular_decomp}
\phi = \phi_+ + \phi_0 + \phi_- 
\end{equation}
where $\phi_+$ and $\phi_-$ extend to holomorphic quadratic differentials on $D^+$ and $D^-$ with simple poles at $z=\infty$ and $z =0$, respectively, and $\phi_0$ is a constant multiple of $\frac{dz^2}{z^2}$.

\begin{lemma}\label{n-int}
Let $\phi$ be a holomorphic quadratic differential on $X_\gamma$. Then
$$\int_{\gamma} \nnn\cdot \phi = \int_{\gamma} \nnn\cdot \phi_0$$
and
$$\left|\int_{\gamma} \nnn\cdot \phi \right| = \ell \|\phi_0\|_\infty$$
\end{lemma}

{\bf Proof:} In the annulus $A_\gamma$ the line integral along $\gamma$ is the line integral over the circle $|z| =1$. 

Note that $\nnn\cdot \phi_\pm$ extends to a 1-form on $D^{\pm}$ so
$$\int_{|z| = 1} \nnn\cdot \phi_{\pm} = 0$$
since by Cauchy's Theorem for any line integral of a closed curve over a holomorphic 1-form on a simply connected region is zero. Therefore
$$\int_{|z| = 1} \nnn\cdot \phi = \int_{|z|=1} \nnn \cdot \phi_+ + \nnn\cdot \phi_0 + \nnn\cdot \phi_- = \int_{|z| = 1} \nnn\cdot \phi_0.$$
We have the covering map $\U\rightarrow A_\ell$ given by $z \rightarrow z^{\frac{2\pi i}{\ell}}$. Thus pushing forward the hyperbolic metric on $\U$ we have that the hyperbolic metric $g_{A_\ell}$ on $A_\ell$ is 
$$g_{A_\ell} = \left(\frac{\ell}{2\pi |z|\cos\left(\frac{\ell}{2\pi}\log|z|\right)}\right)^2 |dz|^2.$$
 As $\phi_0 = \frac{\psi_0}{z^2}dz^2$ we  therefore have 
 $$\|\phi_0\|_\infty = \sup_z \frac{|\phi_0(z)|}{g_{A_\ell}(z)} = \frac{4\pi^2|\psi_0|}{\ell^2}.$$ Thus
$$\int_{|z| = 1} \nnn\cdot \phi_0 = \frac{2\pi i\psi_0}{\ell} \int_{|z| =1} \frac{dz}{z} = \frac{-4\pi^2\psi_0}{\ell}$$ 
and
$$\left|\int_{|z| = 1} \nnn\cdot \phi_0 \right| = \ell\|\phi_0\|_\infty.$$
\eproof

\section{Derivative bounds on univalent maps}
Next we prove Theorem \ref{normbound}. The proof here is considerably more involved. We begin by reducing it to a more explicit statement.
\begin{theorem}\label{main2}
Let $\Phi$ be a quadratic differential on $\U$ and assume that on an $r$-neighborhood of the imaginary axis (in the hyperbolic metric) we have $\|\Phi(z)\| \le K$ with $r< 1/2$ and $K/r < 1/4$. Then there exists a locally univalent map $f\colon \U\to \chat$ such that
\begin{itemize}
\item $Sf = \Phi$;

\item $f(it) \to 0$ as $t\to 0$ and $f(it) \to \infty$ as $t\to \infty$;

\item $f(i) = i$, $|f'(i) - 1| \le \frac{9K}r$ and $|\frac{1}{f'(i)} - 1| \le \frac{9K}r$.
\end{itemize}
\end{theorem}
We now see how this implies Theorem \ref{main}:

{\bf Proof of Theorem \ref{normbound} assuming Theorem \ref{main2}:} Given $ie^t \in \U$ with $t \in \R$ we let $\Phi^t$ be the pull back of $\Phi$ by the isometry $\gamma_t(z) =e^t z$. As $\gamma_t$ preserves the imaginary axis we still have that $\|\Phi^t(z)\| \le K$ for $z$ in an $r$-neighborhood of the axis. Therefore by Theorem \ref{main2} we have a locally univalent map $f^t \colon \U\to \chat$ satisfying the 3 bullets.

We now let $f = f^0$. Since $Sf^t = \Phi^t$ we have that $S(f^t \circ \gamma_{-t}) = \Phi$. Since $f$ and $f^t\circ \gamma_{-t}$ have the same Schwarzian they differ by post-composition with an element of $\pslt$. Since these two maps have the same behavior as $it\to 0$ and $it \to \infty$ this element of $\pslt$ must be of the form $z\mapsto e^\lambda z$. We note that $\nnn$ is invariant under these maps (which includes the maps $\gamma_t$). In particular this implies that
$$\left|\frac{1}{(f^t)'(i)}i - i\right| = \|(f^t)^*(\nnn(i)) - \nnn(i)\|_{g_{\Hp}} = \|f^*(\nnn(f(it))) - \nnn(it)\|_{g_{\Hp}}$$
so the inequality (3) gives
$$\|f^*(\nnn(f(it)) - \nnn(it)\|_{\Hp} \le \frac{9K}r.$$
Similarly if $g_{euc}$ is the Euclidean metric on $\CC^*$ then
$$|(f^t)'(i)i - i| = |(f^t)_*(\nnn(i)) - \nnn(i)| = \|f_*(\nnn(it)) - \nnn(f(it))\|_{g_{euc}}$$
 
\eproof

Here is a brief outline of the proof of Theorem \ref{main2}:
\begin{itemize}
\item Given the quadratic differential $\Phi$ on $\U$ there is an immersion $f_0\colon\U\to \Hs$ such that composition of $f_0$ with the {\em hyperbolic Gauss map} is a map $f\colon \U\to \chat $ with $Sf = \Phi$.

\item The surface $f_0\colon\U\to \Hs$ is the {\em Epstein surface} for $\Phi$. In \cite{epstein-envelopes}, Epstein gives formulas for the metric and shape operator of this surface in terms of the hyperbolic metric on $\U$ and $\Phi$.

\item We will use Epstein's formulas to show that the curve $t\mapsto f_0(ie^t)$ is nearly unit speed and has small curvature. This will imply that $f_0(it)$ limits to distinct points as $t\to 0$ and $t\to \infty$. We then normalize so that these limiting points are $0$ and $\infty$.

\item The proof is then completed by a calculation of the hyperbolic Gauss map using the Minkowski model for hyperbolic space.
\end{itemize}

Before starting the proof of Theorem \ref{main2} we review the necessary facts about {\em Epstein surfaces}. These surfaces are defined for any conformal metric on $\U$ and a holomorphic quadratic differential $\Phi$. Here we will restrict to the hyperbolic metric. 

The {\em projective second fundamental form} for the hyperbolic metric $ g_{\htwo}$ is
$${\rm II} = \Phi + \bar\Phi + g_{\htwo}.$$
The {\em projective shape operator} is given by the formula
$$g_{\htwo}(\hat Bv,w) = {\rm II}(v, w).$$
We then define a {\em dual pair} by
\begin{equation}\label{dualpair1}
g = \frac14\left(\Id + \hat B\right)^* g_{\Hp} \qquad B = \left(\Id + \hat B\right)^{-1}\left(\Id - \hat B\right). 
\end{equation}
By inverting it follows also that
\begin{equation}\label{dualpair2}
g_{\Hp} = \left(\Id +  B\right)^* g \qquad \hat B = \left(\Id +  B\right)^{-1}\left(\Id -  B\right). 
\end{equation}
We will also use two maps on the unit tangent bundle $T^1\Hs$. First we have the projection $\pi\colon T^1\Hs \to \Hs$. We also have the hyperbolic Gauss map $\mathfrak g:T^1\Hs \to \chat$ which takes each unit tangent vector to the limit of the geodesic ray starting from the vector.
\begin{theorem}[Epstein, \cite{epstein-envelopes}]\label{epstein}
Given any holomorphic quadratic differential $\Phi$ on $\U$ there exists a smooth map $\hat f\colon \U\to T^1\Hs$ with
\begin{enumerate}
\item Where $B$ is non-singular the map $f_0 = \pi \circ \hat f$ is smooth, $g = f_0^*g_{\Hs}$ and $B$ is the shape operator for the immersed surface.

\item The map $f = \mathfrak g \circ \hat f$ is locally univalent with $Sf = \Phi$.

\item The eigenvalues of $\hat B$ are $1 \pm 2\|\Phi(z)\|$.
\end{enumerate}
\end{theorem}

\subsection{Geodesic curvature}
The path $\gamma(t) = ie^t$ is a unit speed parameterization of the imaginary axis. We will now begin the computation of the curvature of $\alpha = f_0\circ \gamma$. 

We begin with a few preliminaries. Let $\hat B_0$ be the traceless part of $\hat B$. Then ${\rm II}_0(X,Y) = g_{\Hp}(\hat B_0 X,Y)$ is the traceless part of ${\rm II}$. We then have $\hat B = \Id + \hat B_0$ and ${\rm II}_0 = \Phi + \bar\Phi$. We also need to relate the Riemannian connection $\hat\nabla$ for $g_{\Hp}$ and the Riemannian connection $\nabla$ for $g$. A simple calculation (see \cite[Lemma 5.2]{KS08}) gives  
$$(\Id + \hat B) \nabla_X Y = \hat\nabla_X((\Id +\hat B)Y).$$

\begin{lemma}\label{g_derivative}
Define the holomorphic function $h$ by
$$h = \Phi(\nnn,\nnn).$$
Then
$$\left(\Id + \hat B\right)\nabla_{\dot\gamma}\dot\gamma = 4\Re dh(\dot\gamma)\bar\nnn \quad \mbox{and} \quad \|\nabla_{\dot\gamma}\dot\gamma\|_g = |dh(\dot\gamma)|.$$
\end{lemma}

{\bf Proof:} We work in the complexified tangent bundle of $\U$.  We want to compute $\nabla_{\dot\gamma}{\dot\gamma}$ where $\dot\gamma$ is the tangent vector to the path $\gamma$. As $\left(\Id + \hat B\right)\nabla = \hat\nabla \left(\Id + \hat B\right)$ then
\begin{eqnarray*}
\left(\Id + \hat B\right)\nabla_{\dot\gamma}\dot\gamma &=& \hat\nabla_{\dot\gamma}\left(\Id + \hat B\right) \dot\gamma\\
&=& 2\hat\nabla_{\dot\gamma}\dot\gamma + \hat\nabla_{\dot\gamma}\hat B_0 \dot\gamma\\
&=& \hat\nabla_{\dot\gamma}\hat B_0\dot\gamma
\end{eqnarray*}
since $\hat\nabla_{\dot\gamma}\dot\gamma =0$ as $\gamma$ is a geodesic in $g_{\Hp}$.

To compute this term we note $\dot\gamma = \nnn+ \bar\nnn$. Therefore, since ${\rm II}_0$ is symmetric we have
\begin{eqnarray*}
g_{\Hp}\left(\hat B_0 \dot\gamma, \nnn\right) & =& {\rm II}_0(\nnn +\bar\nnn, \nnn )\\
&=& \Phi(\nnn,\nnn) = h.
\end{eqnarray*}
Using the compatibility of the metric with the Riemannian connection and that $\nnn$ is parallel along $\gamma$ we have
$$ dh(\dot\gamma) = g_{\Hp}\left(\hat\nabla_{\dot\gamma} \left(\hat B_0 \dot\gamma\right), \nnn\right).$$
A similar calculation gives
$$ d\bar h(\dot\gamma) = g_{\Hp}\left(\hat\nabla_{\dot\gamma} \left(\hat B_0 \dot\gamma\right), \bar\nnn\right).$$

Let 
$$v = 2dh(\dot\gamma)\bar\nnn + 2{d\bar h(\dot\gamma)} \nnn.$$
Note that we are extending $g_\Hp$ to be $\CC$-linear on $T\U \otimes \CC$ and therefore
$$g_{\Hp}(\nnn,\nnn) = g_{\Hp}(\bar\nnn,\bar\nnn) = 0 \quad \mbox{and} \quad  g_\Hp(\nnn, \bar\nnn) = 1/2.$$
It follows that
$$g_{\Hp}(v, \nnn) = dh(\dot\gamma) \quad \mbox{and} \quad g_{\Hp}(v, \bar\nnn) = d\bar h(\dot\gamma).$$
As $\nnn$ and $\bar\nnn$ span the tangent space this implies
$$\left(\Id + \hat B\right) \nabla_{\dot\gamma}\dot\gamma= \hat\nabla_{\dot\gamma}\hat B_0(\dot\gamma)=v.$$
Since $\|v\|_{g_{\Hp}} = 2|dh|$ this gives
$$\|\nabla_{\dot\gamma}\dot\gamma\|_{g} = \frac{1}{2} \|\left(\Id+\hat B\right)\nabla_{\dot\gamma}\dot\gamma\|_{g_{\Hp}}=  \frac{1}{2}\|\hat \nabla_{\dot\gamma}\hat B_0\dot\gamma\|_{g_{\Hp}} = |dh(\dot\gamma)|.$$
\eproof

Next we use the Cauchy integral formula to bound $dh$.

\begin{lemma}\label{g_curvature}
Assume that $r \le 1/2$ and $K < 1$. If $\|\Phi(z)\| \le K$ for $z$ in the $r$-neighborhood of $\gamma$ then the geodesic curvature $\kappa_\gamma$ of $\gamma$ in the metric $g$ satisfies
$$\kappa_\gamma \le \frac{5K}{4r(1-K)^2}.$$
\end{lemma}

{\bf Proof:} 
We have
$$\kappa_\gamma = \frac{\|(\nabla_{\dot\gamma}\dot\gamma)^\perp\|_g}{\|\dot\gamma\|_g^2} \le \frac{\|\nabla_{\dot\gamma}\dot\gamma\|_g}{\|\dot\gamma\|_g^2}$$
where $(\nabla_{\dot\gamma}\dot\gamma)^\perp$ is the component of $\nabla_{\dot\gamma}\dot\gamma$ perpendicular to $\dot\gamma$.

To bound $\|\nabla_{\dot\gamma}\dot\gamma\|_g$ we will work in the disk model for $\Hp$ with $\gamma$ the geodesic following the real axis and we will bound the curvature at the origin. With this normalization we have  $\nnn = \frac{1-z^2}{2} \ddz$
and therefore
$$h(z) = \frac{\Phi(z)(1-z^2)^2}{4}.$$
(Here we are not distinguishing between $\Phi$ as quadratic differential and $\Phi$ as holomorphic function.)
At zero, $\dot\gamma = \frac12\ddx$ so by
Lemma \ref{g_derivative} we have
$$\|\nabla_{\dot\gamma}\dot\gamma\|_g = |dh(\dot \gamma)|  = |h_x(0)| /2 = |h_z(0)|/2$$
where the last equality uses that $h$ is holomorphic.

We  now bound  $|h_z(0)|$ using the Cauchy integral formula. Given that $\|\Phi(z)\| \le K$ in the $r$-neighborhood of the origin (in the hyperbolic metric) we have
$$|h(z)| \le \frac14|1-z^2|^2 |\Phi(z)| \le \frac{|1-z^2|^2}{(1-|z|^2)^2} K$$
when $|z| \le R = \tanh(r/2)$. Therefore
\begin{eqnarray*}
|h_z(0)| & =& \frac1{2\pi} \left| \int_{|z|=R} \frac{h(z)}{z^2} dz\right|\\
& \le & \frac{KR}{2\pi R^2(1-R^2)^2} \int_0^{2\pi} \left|1-R^2e^{2i\theta}\right|^2d\theta\\
&=& \frac{K}{2\pi R(1-R^2)^2} \int_0^{2\pi} \left(1-2R^2\cos(2\theta) + R^4\right)d\theta\\
&=& \frac{K(1+R^4)}{R(1-R^2)^2}.
\end{eqnarray*}
If $r\le 1/2$ then, since $\tanh(r/2) \le r/2$, we have $R \le 1/4$. As the derivative of $\tanh^{-1} R$ is $1/(1-R^2)$ when $R \le 1/4$ we have $\tanh^{-1}R \le 16 R/15$ and $\tanh(r/2) \ge 15r/32$. Combining this with our estimate on $|h_z(0)|$ we have
$$|h_z(0)| \le \frac{32}{15 r}\cdot \frac{\left(1+(1/4)^4\right)K}{\left(1-(1/4)^2\right)^2} \le \frac{5K}{2r}$$

and
$$\|\nabla_{\dot\gamma}\dot\gamma\|_g \le \frac{5K}{4r}.$$

We now obtain a lower bound on $\|\dot\gamma\|_g$. For this we recall that by Theorem \ref{epstein} the eigenvalues of $\hat B$ at $z\in \Hp$ are $1\pm2\|\Phi(z)\|$. If $z \in \gamma$ then $\|\Phi(z)\| \le K $ giving 
$$\|\dot\gamma\|_g = \frac12\|\left(\Id +\hat B\right)\dot\gamma\|_{g_{\Hp}} \geq 1-K .$$
Therefore
$$\kappa_\gamma \le \frac{\|\nabla_{\dot\gamma}\dot\gamma\|_g}{\|\dot\gamma\|_g^2} < \frac{5K}{4r(1-K)^2}.$$ \eproof

We can now bound the curvature of $\alpha$ in $\Hs$.
\begin{lemma}\label{curvature_H3}
Assume that $r\le 1/2$ and $K < 1$. Then if $\|Sf(z)\|\le K$ for all $z$ in the $r$-neighborhood of $\gamma$ then the geodesic curvature $\kappa_\alpha$ of $\alpha = f_0\circ \gamma$ in $\Hs$ satisfies
$$\kappa_\alpha \le \frac{3K}{2r(1-K)^2}.$$
\end{lemma}

{\bf Proof:} To bound the geodesic curvature of the curve $\alpha= f_0\circ \gamma$ we again need to bound the covariant derivative $\bar\nabla_{\dot\alpha}\dot\alpha$. This will have a component tangent to the immersed surface which is $(f_0)_*\nabla_{\dot\gamma}\dot\gamma$ plus an orthogonal component whose length is $\|B(\dot\gamma)\|_g$. This last norm is bounded by the product of the maximal eigenvalue of $B$ and $\|\dot\gamma\|_g$. As $B = (\Id+\hat B)^{-1}(\Id-\hat B)$ then eigenvalues of $B$ at a point $f(z)$ in the immersed surface are $-\frac{\|\Phi(z)\|}{\|\Phi(z)\| \pm 1}$. If $z$ is in $\gamma$ we have $\|\Phi(z)\| \le K$ so the maximum eigenvalue is bounded above by $K/(1-K)$. It follows that
$$\kappa_\alpha \le \sqrt{\kappa_\gamma^2 + \left(\frac{K}{1-K}\right)^2 }\le \frac{K}{r(1-K)^2}\sqrt{(5/4)^2+(1/2)^2} \le \frac{3K}{2r(1-K)^2}.$$
\eproof

\subsection{Curves with geodesic curvature $\kappa_g < 1$}
It is well know that curves in $\Hn$ with geodesic curvature $\le \kappa < 1$ are quasi-geodesics. In particular, a bi-infinite path with curvature $\le \kappa$ will limit to distinct endpoints. We need a quantatative version of this statement.
\begin{lemma}\label{quasigeodesic}
Let $\alpha\colon \R\to \Hs$ be smooth, bi-infinite curve with curvature at most $\kappa< 1$. Normalize $\alpha$ in the upper half space model of $\Hs$ so that $\alpha(0) = (0,0,1)$ and $\alpha'(0) = (0,0,\lambda)$ with $\lambda>0$. Then there are distinct points $z_-, z_+\in \chat$ with $\underset{t\to\pm\infty}{\lim} \alpha(t) = z_{\pm}$ with
$$|z_\pm|^{\mp 1} \le \frac1\kappa \left(1-\sqrt{1-\kappa^2}\right)\le \kappa.$$
Let $P_t$ be the hyperbolic plane orthogonal to $\alpha$ at $\alpha(t)$. Then $\del P_t \to z_{\pm}$ at $t \to \pm \infty$.
\end{lemma}

{\bf Proof:} Let $H_\delta$ be a convex region in $\Hs$ bounded by a plane of constant curvature $\delta <1$. If $\alpha'(0)$ is tangent to $\del H_\delta$ and $\delta > \kappa$ then in neighborhood of zero the only intersection of $\alpha$ with $H_\delta$ will be at $\alpha(0)$. Now suppose there is some $t_0$ with, say, $t_0>0$ such that $\alpha(t_0) \in H_\delta$. Then by compactness of the interval $[0,t_0]$ there is a minimum $\epsilon>0$ such that $\alpha$ restricted $[0,t_0]$ is contained in the $\epsilon$-neighborhood of $H_\delta$. This will be a convex region $H_{\delta'}$ bounded by a plane of constant curvature $\delta'>\delta$ and $\alpha$ will be tangent to the boundary of this region. However, all of $\alpha$ restricted to $[0,t_0]$ will be contained in $H_{\delta'}$ contradicting that at the point of tangency the intersection of $\alpha$ and $H_{\delta'}$ will be an isolated point. This implies that $\alpha(0)$ is the unique point where $\alpha$ intersects $H_\delta$ and $\alpha$ is disjoint from the interior of $H_\delta$.

Now take the union of the interior of all regions $H_\delta$ with $\delta>\kappa$ that are tangent to $\alpha'(0)$ and let $\mathcal K$ be its complement.
Then the image of $\alpha$ will be contained in $\mathcal K$.

\begin{figure}[htbp] 
   \centering
   \includegraphics[width=3in]{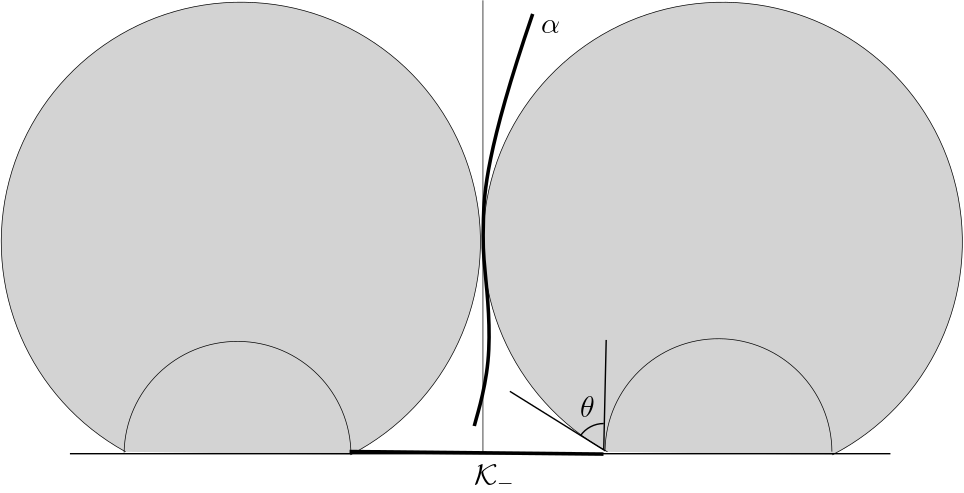} 
   \caption{In this 2-dimensional figure the union of the $H_\delta$ is the shaded region. In $\Hs$ one rotates this region about the vertical axis. The relationship between the angle $\theta$ and a $\kappa$ is given by $\kappa = \sin(\theta)$. One then computes that $\mathcal K_-$ is a disk of radius $\frac{1}{\kappa}\left(1-\sqrt{1-\kappa^2}\right)$.}
   \label{fig:kcurves}
\end{figure}

The intersection of the closure $\mathcal K$ with $\chat = \del \Hs$ will be two regions $\mathcal K_-$ and $\mathcal K_+$ with the accumulation set of $\alpha(t)$ as $t\to\pm \infty$ contained in $\mathcal K_{\pm}$. Then $\mathcal K_-$ is the disk $|z| \le \frac1\kappa \left(1-\sqrt{1-\kappa^2}\right)$ (see Figure \ref{fig:kcurves}) and by symmetry $\mathcal K_+$ is the region with  $1/|z| \le  \frac1\kappa \left(1-\sqrt{1-\kappa^2}\right)$ (including $\infty \in \chat$).

Let $z_-$ be a point in the accumulation set of $\alpha(t)$ as $t \to -\infty$ and let $b_-$ be a Buseman function based at $z_-$. We then observe that the angle between $-\alpha'(t)$ and the gradient $\nabla b_-$ is $\le \theta$.  It is enough to do this calculation when $t=0$. We then let $\mathfrak h$ be the horosphere based at $z_-$ that goes through $\alpha(0)$. Then (along $\mathfrak h$) the gradient $\nabla b_-$ is the inward pointing normal vector field to $\mathfrak h$. Then angle will be greatest when $z_-$ is in $\del \mathcal K_-$ and a direct computation gives that the angle in this case is $\theta$. The bound on the angle implies that $b_-(\alpha(t)) \to \infty$ as $t \in-\infty$ and therefore $\gamma(t) \to z_-$ as $t\to -\infty$.

A similar argument shows that $\alpha(t) \to z_+$ as $t\to\infty$ for some $z_+ \in \mathcal K_+$. Let $\beta$ be the geodesic in $\Hs$ with endpoints $z_-$ and $z_+$. 

Let $z$ be a point in $\del P_t$ and let $\mathfrak h_z$ be the horosphere based at $z$ that goes through $\alpha(t)$ (and is therefore tangent to $\alpha'(t)$). Let $\mathfrak h_{z,\beta}$ be the horosphere based at $z$ that is tangent to $\beta$. We claim that the hyperbolic distance between these two horospsheres is $\le \tanh^{-1}(\kappa)$.  For this we can assume $t = 0$ and let $\mathcal H$ be the convex hull of the regions $\mathcal K_\pm$. Then $P_0 \cap \mathcal H$ is a disk $D$ of radius $\tanh^{-1}(\kappa)$ and $\beta$ will be contained in $\mathcal H$ and intersect $D$. Now let $\mathfrak h_{\pm}$ be the two horospheres based at $z$ that are tangent to $\mathcal H$. The distance between $\mathfrak h_z$ and each of the $\mathfrak h_\pm$ will be $\tanh^{-1}(\kappa)$ and $\mathfrak h_{z,\beta}$ will lie between the $\mathfrak h_\pm$. This implies the distance bound.

We now adjust the picture and so that $z_- = 0$ and $z_+ = \infty$ in the upper half space model for $\Hs$. Assume that $t_i \to -\infty$ and let $z_i$ be points in $\del P_{t_i}$. We'll show that $z_i \to 0$. Assume not. Then, after passing to a subsequence we can assume that $|z_i|$ is bounded below away from zero. Then for each $i$ we apply the isometry $z \mapsto z/|z_i|$ to the points $\alpha(t_i)$ and the horospheres $\mathfrak h_{z_i}$ and $\mathfrak h_{z_i,\beta}$ to get new objects $\bar\alpha(t_i)$, $\bar\hh_{z_i}$ and $\bar\hh_{z_i, \beta}$. This isometry fixes $\beta$ and as $|z_i|$ is bounded below we still have that $\bar\alpha(t_i) \to 0 \in \chat$. Then horospheres $\bar\hh_{z_i, \beta}$ will have Euclidean radius $1$ while the Euclidean radius of the $\bar\hh_{z_i}$ will go to infinity, contradicting that the distance between these horospheres is bounded by $\tanh^{-1}(\kappa)$. Therefore $|z_i| \to 0$, proving the last claim. \eproof

\subsection{Gauss Map}
If $S$ is an immersed surface in $\Hs$ let 
$$n\colon S\to T^1\Hn$$
be the lift to the unit tangent bundle. The Gauss map for $S$ is 
$$\mathfrak g_S\colon S\to \chat$$
with $\mathfrak g_S = \mathfrak g\circ n$.
The following lemma gives the derivative of $\mathfrak g_S$.

\begin{lemma}\label{gauss derivative}
Let $S$ be an oriented, immersed surface in $\Hs$ with $\hat S$ its lift to $T^1\Hs$. Normalize so that $p=(0,1) \in S$ the lift of this point to $\hat S$ is the vector $\ddy$. Let $G\colon T_p S \to T_i\CC$ be the linear map with $G(\ddx) = \ddx$ and $G(\ddt) = \ddy$ and let $\mathfrak g_S\colon S\to \chat$ be the Gauss map. Then
$$(\mathfrak g_S)_*(p) = G \circ (\Id + B)$$
where $B\colon T_p S \to T_p S$ is the shape operator.
\end{lemma}

The proof is a straightforward calculation that is simplest to do in the Minkowski model for $\Hn$.

Let $\langle, \rangle$ be the inner product on $\R^{n,1}$. Then
$$\Hn = \{ x \in \R^{n,1} \ |\  \langle x, x\rangle = -1\}$$
with tangent space
$$T_x\Hn = \{v \in \R^{n,1}\  |\  \langle x, v\rangle = 0\}.$$
Then the restriction of $\langle ,\rangle$ to each $T_x\Hn$ is positive definite and gives a metric of constant sectional curvature equal to $-1$. This is a model for hyperbolic space.

In this model the sphere at infinity for $\Hn$ is the projectivized light cone which we identify with the sphere unit sphere in the plane $x_{n+1} = 1$. In this model there is a very simple formula for the Gauss map.

\begin{lemma}[{Bryant, }\cite{Bryant-Mean}]\label{bryant}
For $x \in \Hn$ and $v \in T^1_x \Hn$ the hyperbolic Gauss map is given by 
$$\mathfrak g(x,v) = \frac{x + v}{\left\langle x+v, (0, \dots, 0,1)\right\rangle}.$$
\end{lemma}

{\bf Proof:} Note that this formula is clear when $x = (0,\dots, 0,1)$. Then general case the follows from equivariance. \eproof

The formula for the Riemannian connection on $\Hn$ is also very simple. The Minkowski connection $\hat\nabla$ on $\R^{3,1}$ is  flat with $\hat \nabla_X Y = X(Y)$. Thus if $\nabla$ is the Riemannian connection on $\Hs$ then by compatibility we have
$$X(Y) = \nabla_XY + \langle X(Y), N\rangle N.$$
where $N$ is the normal to $\Hs$ in $\R^{3,1}$.  We note that $N:\Hs\rightarrow \R^{3,1}$ is given by $N(x) = x$.

From this formula for $\nabla$ we can also calculate the shape operator $B$.
We have
$$BX = \nabla_X(n) = X(n) - \langle X(n), N\rangle N.$$
As $N(p) = p$ then $X(N) = X$ giving
$$\langle X(n), N\rangle = -\langle n, X(N)\rangle = -\langle n, X\rangle = 0$$
so
$$BX = X(n).$$

Note that tangent spaces in $\R^{3,1} \isom \R^{4}$ are canonically identified and
if $p = (0,0,0,1)$ then the linear map $G\colon T_p S\to T_{\mathfrak g_S(p)} \del \Hs$ (from Lemma \ref{gauss derivative}) is the identity map. Therefore, we need to show that
$$(\mathfrak g_S)_*(p) = \Id  + B.$$

With these preliminaries done we can now prove the lemma.

{\bf Proof of Lemma \ref{gauss derivative}:} 
By Lemma \ref{bryant} we have 
$$\mathfrak g_S(x) = \frac{x + n(x)}{\langle x + n(x), (0,0,0,1)\rangle}.$$
Note that for $v \in T_p S$ the derivative of $x \mapsto x + n(x)$ in the direction of $v$ is $v + Bv$ by the calculation above. As $\langle v + Bv, (0,0,0,1)\rangle = 0$ this implies that at $p$ the derivative of the denominator above is zero which gives that
$$(\mathfrak g_S)_*(p) v = v + B v = (\Id + B) v$$
as claimed. \eproof

Using the above, we prove the following bound on the derivative of $f$.
 
\begin{lemma}\label{2Kbound}
Let $f:\U \rightarrow \hat\CC$ be a locally univalent map with Epstein maps $f_0\colon\U\rightarrow \Hs$ and $\hat f_0\colon \U \to T^1\Hs$  normalized as follows:
\begin{enumerate}
\item $f_0(i) = (0,1)$;
\item $\hat f_0(i) = \ddy$;
\item $(f_0)_*(i) \ddy = \lambda \ddt$ with $\lambda = \|\ddy\|_g$.
 \end{enumerate}
Then
$$|f'(i)| = 1 \qquad \mbox{and} \qquad |f'(i) -1| \leq 2\|\Phi(i)\|$$
where $\Phi = Sf$ is the Schwarzian.
\end{lemma}

{\bf Proof:} Let $P \subset T_{(0,1)} \Hs$ be the image of $T_i \U$ under the map $(f_0)_*(i)$. We then have the following sequence of isometries:
$$(T_i\U, g_\Hp) \overset{\Id + B}{\joinrel\relbar\joinrel\relbar\joinrel\relbar\joinrel\longrightarrow} (T_i\U, g) \overset{(f_0)_*(i)}{\relbar\joinrel\relbar\joinrel\longrightarrow} (P, g_\Hs) \overset{G}{\relbar\joinrel\relbar\joinrel\longrightarrow} (T_i\U, g_{\Hp})$$
and by Lemma \ref{gauss derivative} this composition is the map $f_*(i)$. In particular $f_*(i)$ is an isometry of $T_i(U)$ to itself so $|f'(i)| = 1$.

As all the above maps are isometries we can do the computation in any of the metrics. For convenience we will do it in the $g$-metric. Note that $|f'(i) - 1|$ is the distance between the vectors $f'(i)\ddx$ and $\ddx$ in $T_i\U$ or, since $f$ is conformal, the distance between $f'(i) \ddy$ and $\ddy$ with the hyperbolic metric. In the $g$-metric this is the distance between $(\id +B)\ddy$ and $\frac{\ddy}{\lambda}$ (by our normalization (3)). That is
\begin{eqnarray*}
|f'(i) - 1| 
 &= & \left\|(\id+B)\ddy-\frac{\ddy}{\lambda}\right\|_g\\
&= & \left\|(\lambda-1)\frac{\ddy}{\lambda}+B\ddy\right\|_g\\
&\leq& |\lambda-1|\cdot\left\|\frac{\ddy}{\lambda}\right\|_g+\|B\ddy\|_g\\
&=& |\lambda-1|+\|B\ddy\|_g.
\end{eqnarray*}
As $(\id+B)^{-1} = \frac{1}{2}(\id+\hat B)$ and by Theorem \ref{epstein} the eigenvalues of $\hat B$ are $1\pm 2\|Sf(i)\|$ we have 
$$\lambda = \|\ddy\|_g = \frac{1}{2}\|(\id+\hat B)\ddy\|_{g_{\Hp}} \leq 1+\|Sf(i)\|.$$
Therefore
$$|\lambda - 1| \le \|Sf(i)\|.$$
Also $(\Id+\hat B)B = (\Id-\hat B)$ giving
$$\|B\ddy\|_g = \frac{1}{2}\|(\id+\hat B)B\ddy\|_{g_{\Hp}} = \frac{1}{2}\|(\id-\hat B)\ddy\|_{g_{\Hp}}  =  \|Sf(i)\|.$$
Combining these bounds, we obtain the result.
\eproof

\subsection{Proof of Theorem \ref{main2}:}  Let $f^1\colon \U \to \chat$ be the map given by Lemma \ref{2Kbound} and $f^1_0\colon U\to \chat$ the associated Epstein map. By our assumptions $r< 1/2$ and $K/r< 1/4$ so $K< 1/8$. Then Lemma \ref{curvature_H3} gives that
$$\kappa_\alpha < \frac{2K}r<1$$
where $\kappa_\alpha$ is the curvature of $\alpha = f_0^1\circ \gamma$ in $\Hs$.  The normalization (3) in Lemma \ref{2Kbound} and the bound on $\kappa_\alpha$ allows us to apply Lemma \ref{quasigeodesic}. Therefore we have $z_\pm \in \chat$ such that $\underset{t\to \pm\infty}{\lim} f^1_0(it) = z_\pm$ and if $P_t$ are the planes perpendicular to $\alpha$ at $\alpha(t)$ then $\underset{t\to \pm\infty}{\lim} P_t = z_\pm$. Thus  as $f^1(it) \in \partial P_t\cap \chat$ then $\underset{t\to \pm\infty}{\lim}f^1(it) = z_\pm$.

We let $f= m \circ f^1$ where 
$$m(z) = i\cdot \frac{i-z^+}{i-z^-} \cdot \frac{z-z^-}{z-z^+}.$$
Then $m(z^-) = 0, m(z^+) = \infty$ and $m(i) =i$. Therefore $f$ has the desired normalization and we are left to bound the derivative at $i$ which will follow from Lemma \ref{2Kbound} if we can bound the derivative of $m$ at $i$. Computing we have
$$|m'(i)-1| 
=\left| \frac{z^-}{i-z^-} - \frac{i/z^+}{i/z^+ - 1} \right|\le \frac{2\kappa_\alpha}{1 - \kappa_\alpha}.$$
For the reciprocal  we have
$$\left|\frac{1}{m'(i)}-1\right| =\left| \frac{(i-z_-)^2}{i/z_+(z_-/z_+-1)} -\frac{z_-}{i}\right|  \leq \frac{\kappa_\alpha(1+\kappa_\alpha)^2}{1-\kappa_\alpha^2}+ \kappa_\alpha = \frac{2\kappa_\alpha}{1-\kappa_\alpha}.$$
As $\kappa_\alpha < 2K/r$ and by our assumption that  $K/r< 1/4$ then $\kappa_\alpha < 1/2$. Combining these we get that 
$$\left|m'(i)^{\pm 1} -1\right| < \frac{2\kappa_\alpha}{1-\kappa_\alpha} < \frac{8K}r.$$

By Lemma \ref{2Kbound} we have $|(f^1)'(i)| = 1$ and $|(f^1)'(i) - 1| < 2K$ giving
\begin{eqnarray*}
\left|f'(i)^{\pm 1} - 1 \right| &=& \left|((f^1)'(i)m'(i))^{\pm 1} - 1\right| \\
&<&\left |(f^1)'(i)\right|^{\pm 1}\cdot \left|m'(i)^{\pm 1} - 1\right| + \left|(f^1)'(i)^{\mp 1} - 1\right| \\
& < & \frac{8K}{r} +2K < \frac{9K}{r}.
\end{eqnarray*}

\eproof

\section{Application to hyperbolic three-manifolds}
We now prove our application to deformations of hyperbolic 3-manifolds given in Theorem \ref{main_drill} which we now restate.

\medskip

{\bf Theoreom \ref{main_drill}:}
{\em There exists $K,C > 0$ such that the following holds; Let $N$ be an acylindrical hyperbolizable manifold with boundary $S = \partial N$ and $Y\in \T(S)$ with $\|\Phi_Y\|_2 \leq K$.  Then there exists a $\hat Y \in \overline{\T(S)}$ with 
\begin{enumerate}
\item $d_{\rm WP}(Y, \hat Y) \le C\sqrt{\|\Phi_Y\|_2};$

\item $\|\Phi_{\hat Y}\|_\infty \le C \sqrt{\|\Phi_Y\|_2}.$
\end{enumerate}
}

The above result actually holds for {\em relative acylindrical 3-manifolds} as in \cite{wp-paper}. For simplicity we are restricting to the acylindrical setting here. At most points the proof is the same as in \cite{wp-paper}. However, at one key point we will improve the estimate. We will restrict to proving that improved estimate and allow the reader to refer to \cite{wp-paper} for the remainder of the proof.

Let $\cC$ be a collection of essential simple closed curves on $S$ and $\hat N$ the 3-manifold obtained from removing the curves in $\cC$ from the level surface $S\times\{1/2\}$ in a collar neighborhood $S\times [0,1]$ of $N$. Given a $Y\in \T(S)$ there is a unique hyperbolic structure $\hat M_Y$ on $\hat N$ with conformal boundary $Y$. Again there is a projective structure on $Y$ with Schwarzian $\hat\Phi_Y$. Key to our bounds is choosing $\cC$ such that in the complement of the standard collars of $\cC$ we have bounds on the pointwise norm of $\hat\Phi_Y$. It is at this stage that we improve our estimate.

For a complete hyperbolic surface $X$ let $X^{<\epsilon}$ be the set of points of injectivity radius $<\epsilon$. Recall that by the collar lemma (\cite{Buser:collars}) there is an $\epsilon_2 > 0$  such that $X^{<\epsilon_2}$ is a collection of standard collars about simple closed geodesics of length $< 2\epsilon_2$ and cusps. Explicitly we  choose $\epsilon_2 = \sinh^{-1}(1)$, the Margulis constant. If $\gamma$ is closed geodesic of length $<2\epsilon_2$ we let $U_\gamma$ be this collar and for a collection of curves $\cC$ we let $U_{\cC}$ be the union of these standard collars.  We can also choose $\overline\epsilon_2 > 0$ so that for any point in $X^{<\overline\epsilon_2}$ the disk of radius $\inj_X(z)$ is contained in $X^{<\epsilon_2}$. Again, by elementary calculation, we can choose  $\overline\epsilon_2 = \sinh^{-1}(1/\sqrt{3})$.  We denote by $\hat U_\gamma$ and $\hat U_{\cC}$  the corresponding sub-annuli of $U_\gamma$ and $U_{\cC}$. We  note that if $\ell_\gamma(X) < 2\overline\epsilon_2$ then $d(\hat U_\gamma, X-U_\gamma) \geq \overline\epsilon_2$.

We can now state our improved estimate:
\begin{theorem}\label{drill}
There exists constants $C_0, C_1$ such that if $\|\Phi_Y\|_2 \le C_0$ then  exists a collection of curves $\cC$ such that $\ell_\cC(Y) \le 2\|\Phi_Y\|_2$ and
$$\left\|\hat{\Phi}_Y(z)\right\| \le C_1\sqrt{\|\Phi_Y\|}$$ for $z \in Y- \hat U_\cC$.
\end{theorem}

The proof of Theorem \ref{main_drill} follows by combining Theorem \ref{drill} and \cite[Theorem 3.5]{wp-paper}. We now briefly outline this and refer the reader to \cite{wp-paper} for further details.

{\bf Proof of Theorem \ref{main_drill} using Theorem \ref{drill}:} 
In \cite[Theorem 3.5]{wp-paper} we prove that we can find curves $\mathcal C$ such that
\begin{enumerate}
\item $d_{\rm WP}(Y, \hat Y) \le C\sqrt{\ell_{\mathcal C}(Y)};$
\item $\|\Phi_{\hat Y}\|_\infty \le C \sqrt{\ell_{\mathcal C}(Y)}.$
\end{enumerate}
In \cite{wp-paper} we were only able to find $\mathcal C$ such that  $\ell_c(Y) \leq \|\Phi_Y\|_2^{\frac{2}{2n(S)+3}}$ for all $c\in \mathcal C$ which gave the genus dependent bound $\ell_{\mathcal C}(Y) \leq n(S)\|\Phi_Y\|_2^{\frac{2}{2n(S)+3}}$ on the total length of $\mathcal C$. Substituting this bound gave us our original result in \cite{wp-paper}.

To obtain the new bounds in Theorem \ref{main_drill} we apply \cite[Theorem 3.5]{wp-paper} with the improved  bound $\ell_{\mathcal C}(Y) \leq 2\|\Phi_Y\|_2$ from Theorem \ref{drill}.
\eproof

\subsection{$L^2$ and $L^\infty$-norms for quadratic differentials}
In order to prove our estimate, we will first need some results relating the $L^2$ and $L^\infty$-norms for holomorphic quadratic differentials.

We have the following which bounds the pointwise norm in terms of the $L_2$-norm of $\phi$.
\begin{lemma}\label{BW}
Given $\Phi \in Q(X)$ we have:
\begin{itemize}
\item {\bf (Teo, \cite{Teo_curv})} If $z \in X^{\ge \epsilon}$ then 
$$\|\Phi(z)\| \le C(\epsilon) \|\Phi\|_2$$
where $C(x) = \frac{4\pi}{3}(1-\sech^6(x/2))^{-1/2}$.
\item {\bf (Bridgeman-Wu, \cite{BW2021})} If $z \in \hat U_\gamma$ where $\gamma$ is a closed geodesic of length $<2\overline\epsilon_2$ then

$$\|\Phi(z)\| \leq \frac{\left.\|\Phi\right|_{U_\gamma}\|_2}{\sqrt{\inj_X(z)}}.$$
\end{itemize}
\end{lemma}

{\bf Proof:}
The first statement is the main result of \cite{Teo_curv}. The second statement doesn't appear explicitly in \cite{BW2021} but is a simple consequence of it. 
By \cite[Proposition 3.3, part 4]{BW2021} we have for $z \in \hat U_\gamma$
$$\|\Phi(z)\| \leq G(\inj_X(z))\left.\|\Phi\right|_{U_\gamma}\|_2$$
for some explicit function $G$. 
Then by direct computation in \cite[Equation 3.13]{BW2021} we prove that $G(t) \leq 1/\sqrt{t}$ for $t \leq \epsilon_2$. The result then follows.
 \eproof

For a sufficiently short closed geodesic $\gamma$ on $X$ we can lift $\Phi \in Q(X)$ to the annular cover $X_\gamma$.  The covering map is injective on $U_\gamma$  in $X_\gamma$ so 
on $U_\gamma$ we have our annular decomposition  given in \eqref{annular_decomp}
$$\Phi = \Phi^\gamma_- + \Phi^\gamma_0 + \Phi^\gamma_+.$$

We will need the following $L^\infty$-bound in the thin part. In \cite[Lemma 11]{Wolpert:equibounded}) Wolpert proved a bound which is qualitatively the same but we will need the following quantitative version which we derive independently.
\begin{lemma}\label{tail-bound}
Let $z \in \hat U_\gamma$. Then 
$$\|\Phi(z)\| \le \|\Phi^\gamma_0\|_\infty + 2C(\overline\epsilon_2)\left.\|\Phi\right|_{U_\gamma}\|_2$$
\end{lemma}

{\bf Proof:} This is again essentially also contained in \cite{BW2021}. We have on $\hat U_\gamma$ the splitting of $\Phi$ into
$$\Phi = \Phi^\gamma_0 +\Phi^\gamma_+ + \Phi^\gamma_-.$$
Therefore
$$\|\Phi(z)\|\leq \|\Phi^\gamma_0(z)\| + \|\Phi^\gamma_+(z)\|+\|\Phi^\gamma_-(z)\|.$$
We need to bound the norms $\|\Phi^\gamma_0(z)\|$, $\|\Phi^\gamma_+(z)\|$, and $\|\Phi^\gamma_-(z)\|$. By definition $\|\Phi^\gamma_0(z)\| \le \|\Phi^\gamma_0\|_\infty$. By \cite[Proposition 3.3, part 3]{BW2021} 
$$\|\Phi^\gamma_\pm(z)\| \leq C(\overline\epsilon_2)\left.\|\Phi\right|_{U_\gamma}\|_2.$$
Thus it follows that
$$\|\Phi(z)\| \le \|\Phi^\gamma_0\|_\infty + 2C(\overline\epsilon_2)\left.\|\Phi\right|_{U_\gamma}\|_2$$
\eproof

Our next estimate is a statement about holomorphic quadratic differentials that we will use to choose the curves  $\mathcal C$.

\begin{lemma}\label{choose-c}
Let $\Phi \in Q(X)$ with $\|\Phi\|_2 < 1/2$. Let $\mathcal C$ be the collection of simple curves $\gamma$ such that $L_\gamma(X) \leq 2\overline\epsilon_2$ and $\left.\|\Phi\right|_{\hat U_c}\|_\infty \geq \sqrt{\|\Phi\|_2}$. Then 
$$L_{\mathcal C}(X) \leq 2\|\Phi\|_2.$$
 Furthermore for $z \in X-\hat U_{\mathcal C}$ then $\|\Phi(z)\| \leq \sqrt{\|\Phi\|_2}$.
\end{lemma}

{\bf Proof:} For each $\gamma \in \cC$ choose $z_\gamma \in \hat U_\gamma$ with $\|\Phi(z_\gamma)\| \ge \sqrt{\|\Phi\|_2}$. Squaring  the second bound from Lemma \ref{BW} gives the  bound
$$\|\Phi\|_2 \le \|\Phi(z_\gamma)\|^2 \leq \frac{\left\|\left.\Phi\right|_{U_\gamma}\right\|^2_2}{{\inj_X(z_\gamma)}}.$$
Using that $\inj_X(z_\gamma) \ge \ell_\gamma(X)/2$ and summing gives
$$\frac{\|\Phi\|_2}2\sum_{\gamma\in \mathcal C} \ell_\gamma(X) \leq \sum_{\gamma\in \mathcal C}\left.\|\Phi\right|_{ U_\gamma}\|^2_2 \leq \|\Phi\|^2_2$$
which rearranges to give
$$\ell_{\cC}(X) = \sum_{\gamma\in \mathcal C} \ell_\gamma(X) \leq 2\|\Phi\|_2.$$

By the first bound in Lemma \ref{BW} for $z \in X^{\ge \overline\epsilon_2}$ we have as $C(\overline\epsilon_2) \leq 1.1$
$$\|\phi(z)\| \le C(\overline\epsilon_2) \|\Phi\|_2 \le \sqrt{\|\Phi\|_2}$$
where the second inequality uses that $\|\Phi\|_2 \le 1/2$. If $z \in X^{<\overline\epsilon_2}$ but $z \not\in \hat U_\gamma$ for some $\gamma \in \cC$ we have
$$\|\Phi(z) \| \le \sqrt{\|\Phi\|_2}$$
be the definition of $\cC$. These two inequalities give the second bound. \eproof

\subsection{Drilling}
The proof of Theorem \ref{drill} uses the deformation theory of hyperbolic cone-manifolds developed by Hodgson and Kerkhoff in \cite{Hodgson:Kerckhoff:rigidity} and generalized by the second author in \cite{Bromberg:thesis,bromberg}. In particular if $\cC$ is sufficiently short in $M_Y$ (or in $Y$) then there will be a one-parameter family of hyperbolic cone-manifolds $M_t$ where the complement of the singular locus is homeomorphic $\hat N$, the conformal boundary of $M_t$ is $Y$ and the cone angle is $t$.  The $M_t$ will also determine projective structures on $Y$ with Schwarzian quadratic differentials $\Phi_t$ and derivative $\phi_t = \dot\Phi_t$. Then when $t=2\pi$ we have $M_Y= M_{2\pi}$, $\Phi_Y = \Phi_{2\pi}$ and when $t=0$ we have $\hat M_Y = M_0$.

Motivated by the Lemma \ref{choose-c} above, we will fix $\cC$ to be the collection of closed geodesics $\gamma$ on $Y$ with $\ell_\gamma(Y) \leq 2\overline\epsilon_2$ and $\|\left.\Phi_Y\right|_{\hat U_\gamma}\|_\infty \ge \sqrt{\|\Phi_Y\|_2}$. Then by the above  $\ell_{\mathcal C}(Y) \leq 2\|\Phi_Y\|_2$. We first will need to apply the following theorem.

\begin{theorem}[Bridgeman-Bromberg, \cite{BB-L2}]
There exists an $L_0>0$ and $c_{\rm drill}>0$ such that the following holds. Let $M$ be a conformally compact hyperbolic 3-manifold and $\mathcal C$ a collection of simple closed geodesics in $M$ each of length $\le L_0$. Let $ M_t$ be the unique  hyperbolic cone-manifold with cone angle $t \in (0,2\pi]$  about $\mathcal C$ and $\Sigma_t$ the projective structure on the boundary. If $\Phi_t$ is the Schwarzian of the uniformization map for $\Sigma_t$ and $\phi_t = \dot\Phi_t$ then 
$$\|\phi_t\|_2 \le c_{\rm drill} \sqrt{L_\mathcal C(M)}.$$
\end{theorem}

By Bers inequality, $L_{\mathcal C}(M) \leq 2\ell_{\mathcal C}(Y)\leq 4\|\Phi_Y\|_2$. Thus we obtain

\begin{theorem}
If $N$ is an acylindrical hyperbolic manifold with $\|\Phi_Y\|_2 < L_0/4$ then there exists a family of curves $\mathcal C$ such that $\ell_{\mathcal C}(Y) \leq 2\|\Phi_Y\|_2$ and
$$\|\phi_t\|_2 \leq 2c_{drill}\sqrt{\|\Phi_Y\|_2}$$
for $t \in (0,2\pi]$.\end{theorem}

We now promote the $L^2$-bound to a pointwise bound in the complement of the collars of $\mathcal C$.
We first need the following lemma.

\begin{lemma}\label{dbound} There exists a universal constants $D > 0$ such that if $\|\Phi_Y\|_2 < L_0/4$ and $\|\Phi_t(z)\| \le 1/36$ for  all $z \in Y - \hat U_{\cC}$ then  for $z \in Y-\hat U_{\cC}$ we have
$$\|\phi_t(z)\| \le D \sqrt{\|\Phi_Y\|_2}.$$
\end{lemma}

{\bf Proof:} We have the $L^2$-bound $\|\phi_t\|_2 \le 2c_{drill} \sqrt{\|\Phi_Y\|_2}$. We need to convert this to a pointwise bound. 

For $z \in Y^{\ge \overline\epsilon_2}$ we have
$$\|\phi_t(z)\| \le C(\overline\epsilon_2) \|\phi_t\|_2 \le 2c_{drill}C(\overline\epsilon_2) \sqrt{\|\Phi_Y\|_2}.$$

Now assume $\ell_\gamma(Y) \le 2\overline\epsilon_2$ but $\gamma \not\in \cC$. 
For $z \in \hat U_\gamma$ by Lemma \ref{tail-bound} we have
$$\|\phi_t(z)\| \le \left\|\left(\phi_t\right)^\gamma_0\right\|_\infty + 2C(\overline\epsilon_2) \|\phi_t\|_2 \le \left\|\left(\phi_t\right)^\gamma_0\right\|_\infty + 4c_{drill}C(\overline\epsilon_2) \sqrt{\|\Phi_Y\|_2}.$$
We let $D_0 = 4c_{drill}C(\overline\epsilon_2)$. Since the above holds for all $z\in \hat U_\gamma$ we also have
\begin{equation}\label{gamma bound}\|\phi_t\|_\gamma \le  \left\|\left(\phi_t\right)^\gamma_0\right\|_\infty + D_0 \sqrt{\|\Phi_Y\|_2}.
\end{equation}
To finish the proof we need to bound $\left\|\left(\phi_t\right)^\gamma_0\right\|_\infty$ by a constant multiple of $\sqrt{\|\Phi_Y\|_2}$. By definition, for any geodesic $\beta$ with $\ell_\beta(Y) \leq 2\overline\epsilon_2$ then  $d(\hat U_\beta, Y-U_\beta) \geq \overline\epsilon_2$. Therefore $d(\gamma,\hat U_\cC) \ge 2\overline\epsilon_2 > 1/2$ and the $1/2$ neighborhood of $\gamma$ is in $Y-\hat U_{\cC}$. By our assumptions $\|\Phi_t(z)\| \le 1/36$ for $z \in Y - \hat U_{\cC}$ and therefore by Theorem \ref{main} we have
$$\left|\frac{\dot\cL_\gamma(t)}{\ell_\gamma(Y)} + \frac1{\ell_\gamma(Y)}\int_\gamma \nnn\cdot\phi_t\right| \le \frac{9\cdot(1/36)}{1/2} \|\phi_t\|_\gamma = \frac12 \|\phi_t\|_\gamma.$$
Applying Lemma \ref{n-int} to the integral and  rearranging and then applying \eqref{gamma bound} gives
\begin{equation}\label{first bound}\left|\frac{\dot\cL_\gamma(t)}{\ell_\gamma(Y)} \right| \ge \left\|\left(\phi_t\right)^\gamma_0\right\|_\infty - \frac12 \|\phi_t\|_\gamma \ge \frac12\left\|\left(\phi_t\right)^\gamma_0\right\|_\infty - \frac{D_0}{2} \sqrt{\|\Phi_Y\|_2}.\end{equation}

Next we bound the ratio on the left. In \cite{bromberg}, the second author analysed the change in length of geodesics under cone deformations. We now apply a number of results from \cite{bromberg}. By  \cite[Equation 4.6]{bromberg} we have 
\begin{equation}\label{complex bound}
|\dot\cL_\gamma(t)| \le 4 L_{\cC}(t) L_\gamma(t).\end{equation}
We need to bound the two terms on the right. By \cite[Proposition 4.1]{bromberg} we have that $L_\cC(t) \le L_\cC(2\pi)$ for $t\in (0,2\pi]$. When $t=2\pi$ the manifold is non-singular so we can apply the Bers inequality to see that $L_\cC(2\pi) \le 2\ell_\cC(Y)$. Since $|\dot L_\gamma(t)| \le |\dot\cL_\gamma(t)|$ we can combine \eqref{complex bound} with the bound on $L_\cC(t)$ to get
$$\left| \frac{\dot{L}_\gamma(t)}{L_\gamma(t)} \right| \le \left| \frac{\dot{\cL}_\gamma(t)}{L_\gamma(t)}\right| \le 8 \ell_\cC(Y)$$
which integrates to give
$$L_\gamma(t) \le L_\gamma(2\pi) e^{8\pi\ell_\cC(Y)} \le 2\ell_\gamma(Y) e^{16\pi\|\Phi_Y\|_2}.$$
Here the second inequality comes from applying the Bers inequality to $\gamma$ in $M = M_{2\pi}$ and the bound $\ell_\cC(Y) \le 2\|\Phi_Y\|_2$. As we have assumed a universal bound $\|\Phi_Y\|_2 \leq L_0/4$ these two inequalities give
\begin{equation}\label{second bound}\left| \frac{\dot{\cL}_\gamma(t)}{\ell_\gamma(Y)}\right| \le 16 e^{16\pi\|\Phi_Y\|_2} \ell_\cC(Y) \le  D_1\|\phi_Y\|_2.\end{equation}
where $D_1 = 16 e^{4\pi L_0}$ is a universal constant.
Combining \eqref{first bound} and \eqref{second bound} gives
$$ \left\|\left(\phi_t\right)^\gamma_0\right\|_\infty \le (D_0+2D_1) \sqrt{\|\Phi_Y\|_2}.$$
Finally by equation \ref{gamma bound} we have for $z\in \hat U_\gamma$
$$\|\phi_t(z)\| \le  D \sqrt{\|\Phi_Y\|_2}$$
where $D = 2D_0 + 2D_1 = 32e^{4\pi L_0} + 8c_{drill}C(\overline{\epsilon_2})$.
\eproof

Finally we can now prove Theorem \ref{drill} which we  restate.
\medskip

{\bf Theorem \ref{drill} }
{\em There exists constants $C_0, C_1$ such that if $\|\Phi_Y\|_2 \le C_0$ then  exists a collection of curves $\cC$ such that for each $\gamma \in \cC$ we have $\ell_\cC(Y) \le 2\|\Phi_Y\|_2$ and
$$\left\|\hat{\Phi}_Y(z)\right\| \le C_1\sqrt{\|\Phi_Y\|}$$ for $z \in Y- \hat U_\cC$.
}

{\bf Proof:} We let $C_0 = \min\{(1/64)^2, 1/(128\pi D)^2, L_0/4\}$ and $C_1 = 2\pi D+1$.
If $\|\Phi_Y\|_2 < C_0$ then by Lemma \ref{choose-c} for $z \in Y-\hat U_{\mathcal C}$ we have
$$\|\Phi_{Y}(z)\| \leq \sqrt{\|\Phi_Y\|_2} \leq \frac{1}{64}.$$
We now show that  $\|\Phi_t(z)\| \leq 1/36$ for $z \in Y-\hat U_{\mathcal C}$ and $t \in (0,2\pi]$.  If not, then by continuity there is a  $t_0 > 0$ be such that $\|\Phi_t(z)\| = 1/36$ and $\|\Phi_t(z)\| < 1/36$ for $t \in (t_0,2\pi]$. For $t \geq t_0$ as $\|\Phi_Y\|_2 < L_0/4$ by Lemma \ref{dbound} above, we have if $z\in Y-\hat U_{\mathcal C}$ then $\|\phi_t(z)\| \leq D\sqrt{\|\Phi_Y\|_2}$. Therefore integrating we have 
$$\frac{1}{64} \leq \|\Phi_{t_0}(z)-\Phi_Y(z)\| \leq \int_{t_0}^{2\pi} \|\phi_t(z)\| dt < 2\pi D\sqrt{\|\Phi_Y\|_2} \leq \frac{2\pi D}{128\pi D} \leq \frac{1}{64}.$$
This gives our contradiction. Thus we can apply Lemma \ref{dbound} and integrate to get
$$\|\hat\Phi_Y(z)-\Phi_Y(z)\| \leq \int_0^{2\pi} \|\phi_t(z)\| dt \leq 2\pi D\sqrt{\|\Phi_Y\|_2}.$$
Thus
$$\|\hat\Phi_Y(z)\| \leq \|\Phi_Y(z)\|+2\pi D\sqrt{\|\Phi_Y\|_2} \leq (2\pi D +1)\sqrt{\|\Phi_Y\|_2}.$$
\eproof

\bibliography{bib,math}
\bibliographystyle{math}
\end{document}

%% file: macros.tex
\newcommand{\bb}{\mathbb}

\newcommand{\cx}{{\bb C}}

\newcommand{\htwo}{{\bb H}^2}

\newcommand{\ddz}{\frac{\partial}{\partial z}}

\renewcommand{\bold}[1]{\medskip \noindent {\bf \boldmath #1
                        }\nopagebreak[4]}

\newcommand{\del}{\partial}

\newcommand{\isom}{\cong}

\newcommand{\chat}{\widehat{\cx}}

\newcommand{\Hom}{\operatorname{Hom}}

\newcommand{\inj}{\operatorname{inj}}
\newcommand{\id}{\operatorname{id}}

\renewcommand{\Re}{\operatorname{Re}}

\newcommand{\sech}{\operatorname{sech}}

\newtheorem{theorem}{Theorem}[section]

\newtheorem{lemma}[theorem]{Lemma}

\newcommand{\cC}{{\cal C}}

\newcommand{\cL}{{\cal L}}

